\documentclass[11pt]{article}
	
	\newcommand{\blind}{0}
	
    \makeatletter
    \renewcommand\section{\@startsection {section}{1}{\z@}%
                                       {-3.5ex \@plus -1ex \@minus -.2ex}%
                                       {2.3ex \@plus.2ex}%
                                       {\normalfont\fontfamily{phv}\fontsize{16}{19}\bfseries}}
    \renewcommand\subsection{\@startsection{subsection}{2}{\z@}%
                                         {-3.25ex\@plus -1ex \@minus -.2ex}%
                                         {1.5ex \@plus .2ex}%
                                         {\normalfont\fontfamily{phv}\fontsize{14}{17}\bfseries}}
    \renewcommand\subsubsection{\@startsection{subsubsection}{3}{\z@}%
                                        {-3.25ex\@plus -1ex \@minus -.2ex}%
                                         {1.5ex \@plus .2ex}%
                                         {\normalfont\normalsize\fontfamily{phv}\fontsize{14}{17}\selectfont}}
    \makeatother
	\usepackage{amsmath}
	\usepackage{graphicx}
	\usepackage{enumerate}
	\usepackage{xcolor}
	\usepackage{natbib} %comment out if you do not have the package
	\usepackage{url} % not crucial - just used below for the URL
	\usepackage{color}
        \definecolor{darkblue}{rgb}{0.,0.,0.4}
        \usepackage{mathtools, bm, calligra, nameref, tikz}
        \usepackage[font=small,skip=0pt]{caption}
        \usepackage[pdftex,colorlinks=true,linkcolor=black,citecolor=black,urlcolor=black]{hyperref}
	\usepackage[linesnumbered,ruled,lined, noend]{algorithm2e}
        \usepackage{enumitem}[nosep]
        \usepackage{xurl, adjustbox, realboxes,  dsfont, etoolbox,paracol, multirow,   setspace, amssymb, fancyhdr}

        \usepackage{booktabs, makecell,afterpage,  latexsym, multicol}
        \usepackage{rotating}

    \usepackage{xcolor, graphicx, array}
    \usepackage[flushleft]{threeparttable}
    \usepackage{newtxtext}
    \usepackage{float}
    \usepackage{titlesec}
    
    \usepackage{esvect}
    \usepackage{scalerel}
    
    \usepackage{subfigure, epsfig}
    \usepackage{amsmath,amsthm}

\begin{document}
		
			%%%%%%%%%%%%%%%%%%%%%%%%%%%%%%%%%%%%%%%%%%%%%%%%%%%%%%%%%%%%%%%%%%%%%%%%%%%%%%
		\def\spacingset#1{\renewcommand{\baselinestretch}%
			{#1}\small\normalsize} \spacingset{1}
		%%%%%%%%%%%%%%%%%%%%%%%%%%%%%%%%%%%%%%%%%%%%%%%%%%%%%%%%%%%%%%%%%%%%%%%%%%%%%%
		
		\if0\blind
		{
			\title{\bf {An Optimization Framework for Efficient and Sustainable Logistics Operations via  Transportation Mode Optimization and Shipment Consolidation: A Case Study for GE Gas Power}}
			\author{Mustafa C. Camur  $^a$ and Srinivas Bollapragada $^a$ and  Aristotelis  E. Thanos $^b$ \\ and Onur Dulgeroglu $^a$ and  Banu Gemici-Ozkan $^a$\\
			$^a$ General Electric  Research Center, Niskayuna, New York\\
                $^b$ Novelis Inc, Atlanta, GA}
			\date{}
			\maketitle
		} \fi
		
		\if1\blind
		{

            \title{\bf {An Optimization Framework for Efficient and Sustainable Logistics Operations via  Transportation Mode Optimization and Shipment Consolidation: A Case Study for GE Gas Power}}
			\author{Author information is purposely removed for double-blind review}
			
\bigskip
			\bigskip
			\bigskip
			\begin{center}
				{\LARGE\bf An Optimization Framework for Efficient and Sustainable Logistics Operations via  Transportation Mode Optimization and Shipment Consolidation: A Case Study for GE Gas Power}
			\end{center}
			\medskip
		} \fi
		\bigskip
		
	\begin{abstract}
General Electric (GE) Gas Power, a leading manufacturer of gas and steam turbines, manufactures and installs these turbines in power generation plants worldwide. They source components for these turbines from suppliers globally  and transport these components to  manufacturing and assembly locations in the United States using various modes of transportation, including ocean, air, and ground. These transportation options have different lead times and costs. The challenge lies in identifying the most cost-effective solution that meets the assembly requirements, given the high volume of shipments and the complexity of the freight network. To address this challenge,  we  develop a customized, multi-period (dynamic), multi-commodity network flow model and a novel heuristic approach with a rolling time horizon. This model incorporates consolidation and storage options at intermediate nodes, allowing the business to optimize its shipments. %We have successfully implemented and currently utilize this model for managing our global shipments at GE Gas Power.\\
% \textcolor[rgb]{0.00,0.07,1.00}{We strongly encourage authors to address the following three questions in their \textbf{abstract}, preferably following the order shown: (1) Research problem statement: what is the research problem to be addressed? (2) Methods and results: how do the authors address the research problem and what are the main results? (3) Insights and implications:  What have the authors learned (as opposed to what they did, which is covered in point (2)) from conducting this research? What is the knowledge gained and why does it matter? The abstract should be written in \textbf{a single paragraph}.}.\\
% We thank you for your attention to these details.
	\end{abstract}
	\noindent%

	%\newpage
	\spacingset{1.5} % DON'T change the spacing!

\section{Introduction} \label{Introduction}

General Electric (GE) Gas Power, a business of GE company, is the world's leading manufacturer of gas and steam turbines and provider of related  services. It has a large installed base of gas turbines, steam turbines, and generators, collectively  generating  approximately one-third of the world's electricity  \citep{Exclusiv95}.    As the world shifts towards a lower carbon future, GE Gas Power's mission is to lead the energy transition by pioneering technologies such as hydrogen fuel usage, carbon capture, and sequestration to limit climate change. %They also leverage natural gas to support renewable technologies to provide sustainable, affordable, and reliable electricity to more people worldwide. 

\begin{figure}[!htbp]
\centering
     \caption{Visualization of the product order flows via different transportation modes}
     \includegraphics[page=1,width=1 \linewidth,keepaspectratio=true]{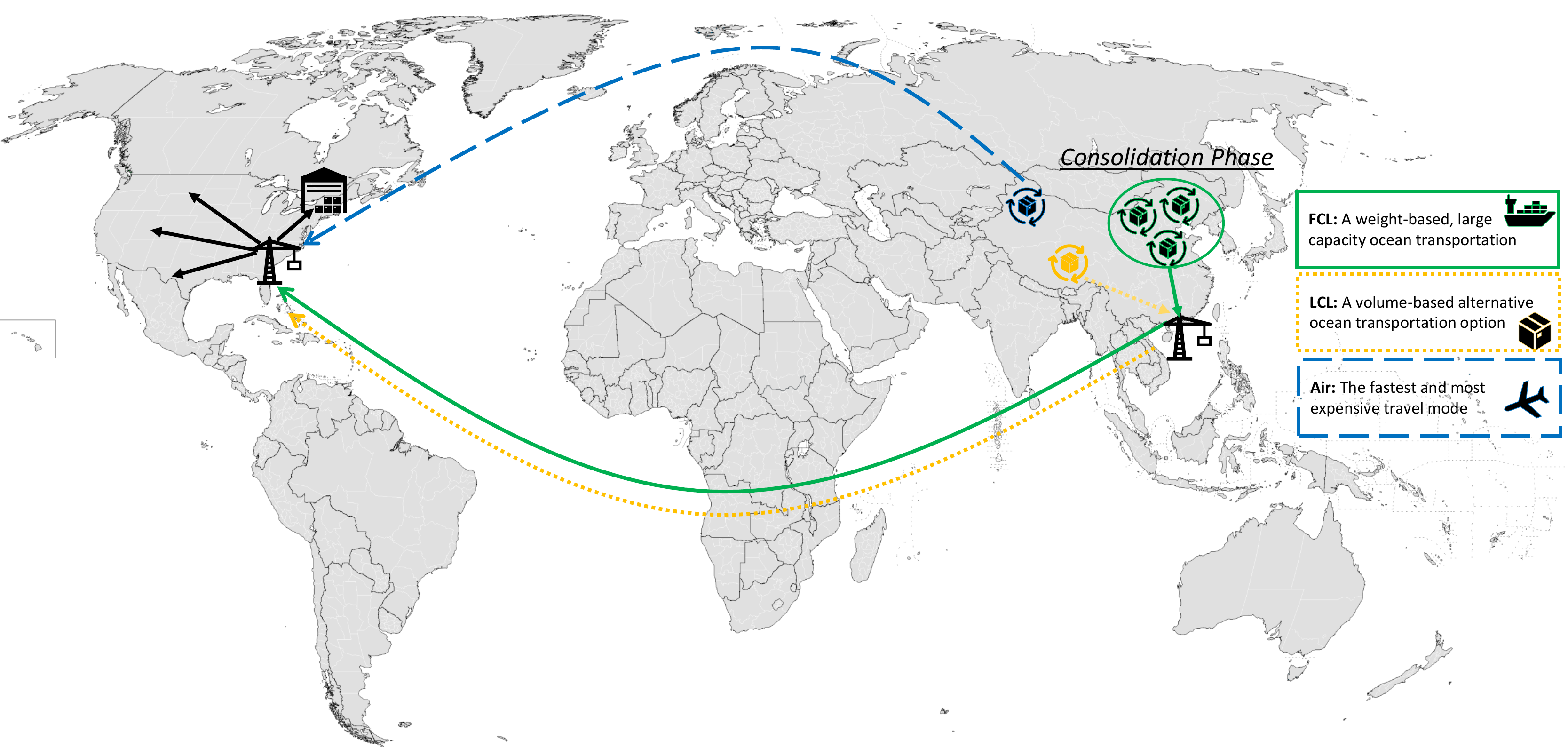}   \label{Figure1}
\end{figure}

GE Gas Power manufactures gas and steam turbines at its plant in Greenville, South Carolina. The turbine components come from suppliers worldwide and are transported to the manufacturing facilities using various modes of transportation, including ground, ocean, air, or a combination of these methods \citep{camur2023enhancing}. Three transportation options exist for moving parts from overseas suppliers to the manufacturing plants in the United States (U.S.), each with different transportation costs and lead times. We describe these options below and illustrate them in Figure~\ref{Figure1}.
\begin{enumerate}
    \item  Full Container Load (FCL) Shipment: In this option, GE books an entire container to transport the parts. The capacity of a container ranges between 20 and 27 tons, and the delivery time varies between 5 and 7 weeks. Reserving an FCL container needs to be done a few weeks in advance, and a fixed cost for the container is incurred regardless of the amount of cargo loaded. Therefore, when there are enough parts available to sufficiently fill the container, the FCL shipment option has the lowest cost compared to other options. Additionally, since multiple parts are consolidated and shipped in the same container, it is easier to track them during the transportation phase.
    
    \item 	Less than Container Load (LCL) Shipment: The cost of shipping items under this option depends on their volume, resulting in higher costs for larger units. The transportation lead time for LCL shipments is also between five and seven weeks, similar to the FCL option. However, unlike FCL, it is not necessary to reserve space in the container several weeks in advance. The LCL option is preferred for ad-hoc and lower-volume shipments. To utilize LCL shipping, the company needs to enter into an agreement with a third-party carrier, with the fee depending on the shipment's volume. It's important to note that in the LCL option, no single company owns all the contents in the container. Instead, products from different companies are consolidated within the same container. This option is preferred when there are insufficient parts available to fill an entire container, and there is enough time for transporting the parts.
    \item 	Air Shipment: It is the fastest method of transporting a product and is generally the most expensive option in almost all situations. Therefore, air shipment is primarily preferred when a product has a tight deadline or has already passed it. It is also chosen when ocean transportation options such as FCL and LCL are not feasible for transporting the product.
\end{enumerate}

% \begin{figure}[!htbp]
% \centering
%      \caption{ \textbf{Full Container Load} where all of the shipped packages belong to the same company vs. \textbf{Less Than Container Load} where the logistics carrier consolidates shipments from various companies into the same container \citep{LessThan16} }
%      \includegraphics[width=0.75\linewidth,keepaspectratio=true]{FCLvsLCL.png}   \label{Fig1}
% \end{figure}

Furthermore, there is a preferred time window for delivering parts to the destination location. It is important to avoid delivering the shipments too early to prevent inventory carrying costs and potential unavailability of unloading equipment (e.g., cranes). In addition, they cannot be late  as this could lead to delays in the manufacturing and assembly operations at the destination, resulting in significant costs  \citep{nguyen2021towards}. The choice of the transportation option not only affects the transportation cost but also determines the arrival time at the destination location.

The importance of using consolidation operations  extends beyond cost savings. While shipping a product by container results in lower carbon emissions compared to air shipping, container shipments still have a significant impact on carbon emissions \citep{pettit2018revisiting}, which are the primary drivers of climate change \citep{lee2021scale}. Underutilization of containers leads to a higher carbon footprint, contributing negatively to global warming \citep{ma2019reducing}. By consolidating shipments, more sustainable operations can be achieved, aligning with GE's clean operation goals. The supply chain and logistics operations of GE Gas Power's business, like those of many other companies, have been impacted by the global pandemic and conflicts worldwide \citep{camur2023integrated}. To improve the transportation operations of these parts, the GE Research team collaborated with the Gas Power team to identify areas in need of enhancement.
%Previously at GE Research Center, researchers have investigated predictive modeling in other contexts of optimization including manufacturing \citep{carter2019open, felix2022situ} and process optimization \citep{salasoo2020system, ashenafi2022reinforcement} teaming up with different GE businesses. 

FCL consolidation operations have only recently been  started at GE Gas Power. In the past, the business solely relied on LCL and air shipment options.  In addition, the business did not have a decision support system in place to effectively choose among the three transportation options for parts.  Consequently,  some FCL units were shipped without fully utilizing their capacity, while other FCL opportunities were completely overlooked. Additionally, a significant number of products were air shipped due to time constraints and deadlines, leading to elevated transportation costs. The underutilization of FCL and LCL containers resulted in increased carbon emissions.

In this study, we develop an algorithm and a software tool to  determine the optimal shipment options for all the parts to be transported.  We formulate the parts shipment problem as an integer linear program inspired by the dynamic multi-commodity network flow problem and solve it using the state-of-the-art Gurobi solver. While the commercial solver can solve smaller problems to optimality in a reasonable amount of time, the large-scale real-world problems require more  time and return larger optimality gaps. We, therefore, develop a novel heuristic algorithm that is based on the well-known Knapsack Problem to solve the problem in a few minutes (see \cite{camur2023stochastic} for a similar approach).

It is crucial to note that our heuristic approach primarily  focus on identifying the optimal timing for  FCL booking decision, especially in the context of our long-term planning. As  new product orders come in and / or some product orders get canceled, the predetermined routing and placement decisions might not always align with the business's needs. Hence, our heuristic approach helps the business with determining when, where and how many FCL units to book. 

The remainder of the paper is organized as follows. We first present a literature review in Section \ref{Literature Review}. In Section \ref{Mathematical Optimization Model}, we describe our integer programming model. The heuristic approach to solve the problem is discussed in Section \ref{A Knapsack-based Heuristic Approach}. We conduct an extensive computational study to compare different methods and analyze the results in Section \ref{Computational Experiments}. Finally, in Section \ref{Conclusion}, we summarize our work and discuss potential  future research directions.

\section{Literature Review}
\label{Literature Review}
In this section, we review the existing literature that is primarily  related to our work and discuss how our study differs from previous works while contributing to the literature. Our problem is mainly associated  with the multi-commodity network flow problem and its variants, as well as the studies on shipment consolidation.

\subsection{Multi-Commodity Network Flow}
\label{Multi-commodity Network Flow Problem}
The multi-commodity network flow problem (MCNFP) finds a place in the literature starting from the late 1950s \citep{ford1958suggested}. 
%\cite{ford1958suggested} initially proposed the problem, which was later generalized by \cite{hu1963multi}. %Since then, researchers have worked on the problem in various application areas, including network interdiction \citep{lim2007algorithms}, fashion retailing \citep{zhang2019multi}, liner shipping \citep{trivella2021multi}, and disaster response \citep{chang2022multi}. 
\cite{borradaile2017multiple} analyzed the MCNFP on directed planar graphs and presented a solution algorithm running in pseudo-linear time. \cite{ouguz2018multicommodity} studied multi-commodity flows on a restricted continuous location problem and designed a decomposition-based exact solution technique. They considered both deterministic and stochastic versions to evaluate the effectiveness of their algorithm. \cite{xu2018integrated} formulated an integrated train timetabling and locomotive assignment problem for a multi-state railway network as a minimum cost MCNFP. The authors employed a Lagrangian relaxation-based heuristic to solve their problem at scale, providing a practical solution approach.

 \cite{castro2020mdd} focusd on the multi-commodity Travelling Salesman Problem and  proposed a solution approach using multi-valued decision diagrams. The problem addresses the delivery of multiple products to their respective 
 destinations while ensuring that the capacity of each vehicle is not exceeded. In a different application domain,  \cite{wang2022new} proposed a compact MCNFP model to enhance the operational capabilities of an airport.  Their model aims to optimize both robustness and taxiing times, contributing to the efficiency of airport operations. \cite{salimifard2022multicommodity} also presented a review paper on the multi-commodity network flow problem, summarizing the recent advancements in the field. In our study, we specifically focus on the time-expanded version of the MCNFP and incorporate business-related constraints into our model. By addressing this specific variant and considering real-world constraints, we aim to contribute to the existing literature on MCNFP.

%Furthermore, several survey and review papers have been published on the topic of MCNFP. For instance, \cite{kennington1978survey} provided a comprehensive survey of the problem, discussing various solution approaches and applications. 

\subsection{Multi-Period Multi-Commodity Network Flow} \label{Multi-Period Multi-Commodity Network Flow}

The multi-period, multi-commodity network flow problem is also well studied \citep{aronson1989survey, camur2021optimizing}. It is also known as the dynamic MCNFP (DMCNFP) and MCNFP in time-varying networks. 
%Since the literature on this topic is remarkably vast, we will only focus on a few studies related to our work in this section.   
\cite{zhang2017time} discussed bicycle sharing systems modeled as time-space network flow models. The authors employed forecasting methods and incorporated vehicle routing decisions into their model. \cite{trivella2021multi} modeled an MCNFP formulation where transit times are treated as soft constraints, penalizing shipment delays in the objective function. They used a path-based formulation for which a column generation approach is proposed as an exact solution technique. However, neither of these studies allows for inventory/storage at the intermediate nodes, which is a crucial aspect in our problem setting due to parts consolidation operations.

%\cite{bertsimas2000traffic} analyzed air traffic control under changing weather conditions, aiming to minimize the costs of delays. They utilized Lagrangian relaxation to aggregate certain flows and then employed a randomized rounding heuristic to create a set of flight paths. \cite{glockner2000dynamic} studied the DMCNFP in a stochastic environment, where arc capacities are treated as random variables, and proposed various decomposition schemes. The authors designed a polynomial time algorithm for the subproblems generated after decomposition.

In a recent study  on MCNFP in time varying networks, \cite{khodayifar2021minimum} proposed a formulation that allows for inventory at intermediate nodes, along with a decomposition principle. However, this study did not consider consolidation operations and did not analyze the impact of different transportation modes. Yet, multi-modal transportation is necessary and often unavoidable in the majority of international shipment operations. This introduces additional challenges in optimization models and solution techniques \citep{archetti2022optimization}. 

In most of the literature on both MCNFP and DMCNFP, flow variables are treated as continuous variables but we cannot split our products into fractional parts. Hence, we use only integer variables in our formulation, significantly increasing the complexity of the problem. \cite{gendron2017reformulations} introduced reformulations for MCNPF via discretization, thereby exploiting integral flows. The authors also introduced valid inequalities, cut-set inequalities, and flow disaggregation methods to strengthen their formulation. However, this technique has not been extended to DMCNFP formulations.

\subsection{Multi-Commodity Network Flow and Network Design}
\label{Multi-Commodity Network Flow and Network Design}

Both MCNFP and DMCNFP have extensive applications in network design problems as well \citep{belieres2020benders, hewitt2022flexible}.  Even though our specific problem setting does not include a full network design component, the FCL booking decisions have an impact on whether the corresponding arc can be used for transportation operation. Thus, we will cover some relevant multi-commodity network design (MCND) studies in this context.

%\cite{karsten2015time} studied the liner shipping network design problem aiming to determine which cargoes should be transported via which routes. The authors incorporated transit times into their models to capture a more realistic scenario. Building upon this work,  \cite{karsten2017competitive}analyzed a similar network design problem with different service levels in cargo shipments. They employed an efficient meta-heuristic algorithm as a solution method to tackle the problem in real-world liner shipping networks. 
\cite{kazemi2021new} introduced new commodity representations for the MCND problem, allowing partial aggregation of commodities. The authors designed two new optimization models based on partial aggregations and conducted a polyhedral analysis to compare the LP relaxations of the new and existing formulations. In a different setting, 
 \cite{kazemzadeh2022node} analyzed the multi-commodity flow in a fixed-charge network design problem and proposed a node-based Lagrangian relaxation approach embedded into a metaheuristic. Computational experiments demonstrated that the proposed solution method competes favorably with existing heuristics in the literature.

\cite{li2019post} developed a dynamic, bilevel MCND model to solve a post-disaster road network repair problem, where repair crews and relief materials are treated as commodities. They utilized a tailored-genetic algorithm to solve their model. 
\cite{ghasemi2019uncertain} studied an evacuation planning problem, which is formulated as a dynamic MCND problem with multi-objectives under a stochastic environment. The goal is to use multiple vehicles to transport injured people while providing them with relief commodities. 

\cite{kirschstein2019multi} discussed dynamic MCND for a lot-sizing and storage selection problem, incorporating both supplier and storage selection in the planning phase. They employed a kernel-search-based heuristic to solve large-scale test instances.   \cite{fragkos2021decomposition} developed new decomposition algorithms for the dynamic MCND problem, where the design aspect is deemed as network expansion. The authors proposed Lagrangian relaxation supported by a heuristic and a modified Benders decomposition for the capacitated and uncapacitated versions of the problem, respectively. However, none of these studies took the consolidation operations and the resulting cost-saving opportunities into consideration. Therefore, in the next section, we will examine studies that specifically cover inbound shipment consolidation.

\subsection{Shipment Consolidation}
\label{Shipment Consolidation}
Shipment consolidation is a critical transportation strategy employed by companies to reduce transportation costs. We refer the reader to the doctoral dissertation by \cite{fang2006inbound} where the author discusses how consolidation may help transshipment companies working between China and the U.S. with significant cost savings  (see also \citep{tyan2003evaluation}). Moreover,   \cite{wei2021shipping} examined the economy of scale opportunities associated with consolidation, specifically FCL and LCL, and explore the trade-off between consolidation and air transportation when faced with tight deadlines.

%One of the  earliest papers  on inbound consolidation is  by \cite{gupta1987inbound}, where the authors discuss the benefits of consolidation in logistics operations, specifically in terms of reducing shipment costs. 
\cite{popken1994algorithm} proposed a nonlinear MCNFP model containing freight consolidation for manufacturing companies to decrease transportation cost while maintaining a high level of service performance. The author developed an iterative approach that utilizes a linearization method to identify local optimal solutions and a heuristic method to improve local solutions.  
\cite{zhang2021multi} designed an MCNFP model including package consolidation in the online retailing business and devised a logic-based Benders decomposition method to solve the problem. In this study, the authors examined the trade-off between consolidation and splitting. Yet, none of these studies considered multi-period planning or addressed the business-related constraints that we are bound by. %In addition, there is a plethora of studies including consolidation operations in ground transportation (e.g., less than truck load); however, since our challenge primarily focuses on ocean consolidation, we will not delve further into this topic.

 \cite{hanbazazah2019freight} studied the DMCNFP problem with both FCL and LCL consolidation for a third-party logistics company. They proposed a mixed-integer linear programming model and employed a three-phase decomposition technique to solve it. While this study is the most similar to our problem, there are still significant differences. Firstly, the authors did not consider alternative transportation methods other than FCL/LCL, and every shipment was forced to go through the gateways for consolidation. In contrast, we allow for air shipments between airports, especially in cases where ocean shipments may not meet deadlines, necessitating more expedited air shipment. Additionally, the model in the mentioned study did not enforce an advance booking constraint or limit the number of containers that can be booked. However, in reality, FCL units must be booked several weeks in advance, and there may be restrictions on the number of containers that can be booked, particularly during periods of strained global supply chain networks. Finally, the authors did not limit storage duration at transit (i.e., intermediate) nodes, allowing a product to wait at a gateway for as long as needed for consolidation. Yet, in many ports, this is not realistic due to inventory and tracking-related considerations. Therefore, we believe that our work adds realism to the existing models in the literature, enabling practical implementation in real-world scenarios.

\section{Problem Formulation} \label{Mathematical Optimization Model}

In this section, we introduce an optimization model to solve the inbound logistics problem under time window constraints for part-delivery times, where storage of products is permitted at intermediate nodes. The sets, parameters, and variables are presented in Tables \ref{Table:1}, \ref{Table:2}, and \ref{Table:3}, respectively. Following this, we present our formulation with a detailed explanation of each constraint.

\begin{table}[ht]
\centering
\caption{Set definitions}
\setlength\tabcolsep{4pt} 
\scalebox{1}{
\begin{tabular}{ll} 
 \hline
 Set & Definition  \\ 
 \hline
$V$ & locations (i.e., nodes) \\
$S$ & source nodes where $S \subset  V$ \\
$I$ & in-transit nodes where $I \subset V$ \\
$D$ & demand nodes where $D \subset  V$ \\
$E$ & connections between locations (i.e., edges) \\
$G = (V,E) $ & transportation network \\
$M^E$ & transportation modes  associated with each edge \\
$M^A$ & air transportation with corresponding  edges \\
$M^L$ & LCL transportation with corresponding  edges \\
$M^F$ & FCL transportation  with corresponding  edges \\
$M^G$ & ground transportation  with corresponding  edges \\
$P$ & product orders  \\
$T$ & time horizon in weeks\\
\hline
\end{tabular}}
\label{Table:1} 
\end{table}

\begin{table}[ht]
\centering
\caption{Parameter definitions}
\setlength\tabcolsep{4pt} 
\scalebox{1}{
\begin{tabular}{ll} 
 \hline
 Parameter & Definition  \\ 
 \hline
 $c_{ijmp}$ & transportation cost of product order $p$ from location $i$ to location $j$ via mode $m$ \\
$f_{ij}$ & fixed cost of using FCL \\
% $f_{ij}^{L}$ & fixed cost of using LCL \\
$\tau_{ijm}$ & time taken to travel from location $i$ to location $j$ via mode $m$ \\
$\sigma_{ijm}$  & capacity of edge $(i,j)$ when using  mode $m$ \\
$d_{pi}$ & whether product order  $p$ exists at demand node $i$\\
$s_{pi}$ & whether  product order $p$ comes from supply node $i$\\
$a_{pi}$ & availability time of product order $p$ at supply node $i$\\
$\xi_{p}$ & early deadline of delivering product order $p$\\
$\Xi_{p}$ & late deadline of delivering product order $p$ \\
$w_p$ & weight  of product order $p$ \\ 
$\rho$ & number of weeks that a product order can remain idle at an in-transit node \\
$\upsilon$ & number of weeks that an FCL unit is required to be book in advance\\
$\lambda $ & number of FCL units allowed to be booked at a port per week\\
\hline
\end{tabular}}
\label{Table:2} 
\end{table}

\begin{table}[ht]
\centering
\caption{Variable definitions}
\setlength\tabcolsep{4pt} 
\begin{tabular}{ll} 
 \hline
 Variable & Definition  \\ 
 \hline
$x_{ijmpt} $ & whether product order $p$ is sent from location $i$ to location $j$ via  mode $m$  at time $t$\\ 
$z_{ijmt}$ & whether an FCL unit $m \in M^F$ is used going from location $i$ to location $j$ at time $t$ \\
$e_{ip}$ & arrival time of product order $p$ at demand location $i$\\
\hline
\end{tabular}
\label{Table:3} 
\end{table}

\allowdisplaybreaks
\begin{subequations} 
\begin{align}
 \label{obj}  
 \textbf{\hypertarget{IP}{IP}:} \nonumber \\ ~\min &~ \sum_{ (i,j) \in E} \sum_{p \in P}  \sum_{t \in T} \sum_{m \in M} c_{ijmp}  x_{ijmpt} +  \sum_{ (i,j) \in E} \sum_{m \in M^{F}}   \sum_{t \in T}  f_{ij} z_{ijmt}
  \\
\label{capacity}
s.t. &~ \sum_{p \in P} w_p x_{ijmpt} \leq \sigma_{ijm}, \quad  (i,j)  \in E, m \in M^E,  t \in T \\
\label{supply}
&~    \sum_{ j: (i,j) \in E} \sum_{m \in M^E } \sum_{t\in T} x_{ijmpt} \leq s_{pi}, \quad   i \in S,  p \in P \\
\label{in-transit-1}
&~  \sum_{ j : (j,i) \in E} \sum_{m \in M^E } \sum_{t - \tau_{jim} -\rho \leq t^{\prime} \leq t - \tau_{jim} } x_{jimpt^{\prime}}  
\geq  \sum_{ k: (i,k) \in E} \sum_{m \in M^E } x_{ikmpt}, \quad i \in I,  p \in P,    t \in T  \\
\label{in-transit-2}
&~ \sum_{t-\rho \leq t^{\prime} \leq t} \sum_{ j : (j,i) \in E} \sum_{m \in M^E }    x_{ijmpt^{\prime} }  
\geq  \sum_{ k: (i,k) \in E} \sum_{m \in M^E } x_{kimp(t - \tau_{jim} -\rho)}, \quad i \in I,  p \in P,    t \in T  \\
\label{in-transit-3}
&~  \sum_{ j : (j,i) \in E} \sum_{m \in M^E } \sum_{t \in T }   x_{jimpt}  
=  \sum_{t \in T }  \sum_{ k: (i,k) \in E} \sum_{m \in M^E } x_{ikmpt}, \quad i \in I,  p \in P  \\
 \label{demand}
&~    \sum_{j: (j,i) \in E} \sum_{m \in M^E } \sum_{t\in T} x_{jimpt} = d_{pi}, \quad   i \in D,  p \in P \\
% \label{linking}
% &~ x_{ijmpt} \leq  y_{ijmpt}, \quad   (i,j) \in E,  p \in P, m \in M^E,  t \in T  \\
\label{FCLusage}
&~ x_{ijmpt} \leq z_{ijmt} , \quad  m \in M^F, (i,j) \in E,  p \in P,    t \in T  \\
\label{restrictFCL}
&~ z_{ijmt} =0 , \quad   m \in M^F, (i,j) \in E,    t < \upsilon : t \in T  \\
\label{numberofFCLsperPort}
&~ \sum_{j : (i,j,) \in E } \sum_{m \in M^F}z_{ijmt} \leq \lambda, \quad   i \in  I, t \in T  \\
\label{MCF5}
&~  \sum_{j : (i,j) \in E} \sum_{m \in M^E} \sum_{t \in T: t < a_{pi}} x_{ijmpt}  = 0, \quad  i \in S,  p \in P \\
\label{traveltime1}
&~ (t  +  \tau_{ijm}) x_{ijmpt} \leq e_{jp}, \quad   (i,j) \in E : j \in D,  p \in P,  m \in M^E, t \in T  \\
\label{traveltime2}
&~( t  +  \tau_{ijm}) +  M (1-x_{ijmpt}) \geq e_{jp}, \quad   (i,j) \in E : j \in D,  p \in P,  m \in M^E, t \in T  \\
\label{deadlines}
&~ \xi_{p}   \leq e_{ip} \leq \Xi_{p}  , \quad   i \in D, p \in P \\
\label{variables}
&~ x_{ijmpt} \in \{0,1\}, z_{ijmt} \in \{0,1\}, e_{ip} \in \mathbb{Z}_{+} &
% \epsilon_{ip} \in \mathbb{R}_{+}, \ell_{ip} \in \mathbb{R}_{+}   
\end{align}
\end{subequations}

The objective function \eqref{obj} aims to minimize i) the total transportation cost, and ii) the total cost of the FCL units used. For the rest of the paper, each product order is referred to as a product for simplicity. We ensure that no edge capacity is violated through Constraints \eqref{capacity}. Constraints \eqref{supply} and \eqref{demand} are similar to the generic multi-commodity flow conservation constraints under time windows for supply and demand nodes, respectively. Constraints \eqref{in-transit-1}, \eqref{in-transit-2}, and \eqref{in-transit-3} are the flow balance constraints for in-transit nodes where products are permitted to be stored if necessary. This is a realistic assumption as some products might arrive at an in-transit node early and wait to be consolidated into an FCL unit. Note that a product cannot remain at an in-transit node for more than $\rho$ weeks.

 Through Constraints \eqref{FCLusage}, we ensure that even if one product is sent via an FCL unit, it triggers a fixed cost in the objective function. No FCL unit can be utilized in the first $\upsilon$ weeks, as ensured by Constraints \eqref{restrictFCL}. This is because it takes $\upsilon$ weeks to book an FCL unit in reality, based on the agreement signed between the business owner and the shipment company. In addition, Constraints \eqref{numberofFCLsperPort} ensure that the number of FCL units booked at a port is limited.

Constraints \eqref{MCF5} guarantee that the arrival time of each product at a source node is set as the product's availability time. Both Constraints \eqref{traveltime1} and \eqref{traveltime2} compute the arrival time of each product at its destination. It can be observed that deadline constraints could be captured through variables $x_{ijmpt}$ by conditioning on the time index $t$. However, as we will discuss in the next section, there could be scenarios where deadlines might be relaxed for some business cases. Hence, we define variables $e$ for potential modeling changes. Lastly, Constraints \eqref{variables} place the bound restrictions on our variables.

\subsection{Modelling Improvements} \label{Model Improvements}
Introducing both new variables and constraints, including valid inequalities and constraint tightening, are shown positive impacts in several  studies \citep{camur2021large, luo2022repeated}. This application requires that all products be delivered within their respective earliest and latest deadlines. Constraints \eqref{deadlines} that enforce these deadlines can sometimes lead to infeasibility if an order is not available for shipment in sufficient time before its latest deadline. We handle such orders outside the model.  They are typically shipped by the business through air shipment or other means.

It sometimes might be acceptable to miss deadlines while incurring penalty costs. In such situations, an alternative is to relax these constraints and use penalty variables to penalize failure to meet the constraints. Let $\epsilon_{p}$ and $\ell_{p}$ be the variables representing the penalty costs for product $p$ arriving early and late at its demand node, respectively. We can rewrite Constraints \eqref{deadlines} as:

\begin{equation} \label{relaxedDeadline}
    \xi_{p}  - \epsilon_{p} \leq e_{ip} \leq \Xi_{p} + \ell_{p}   , \quad   i \in D, p \in P
\end{equation}

Constraints \eqref{relaxedDeadline} indicate that if no feasible solution exists for a product with a given input, either early or late deadlines can be violated. Additionally, we may add $\sum_{p \in P} \sum_{i \in D} (\epsilon_{p}+ \ell_{p})$ to the objective function to ensure that the model does not arbitrarily relax the deadline constraints. %This modification turns our problem into a multi-criteria decision-making model.

In other applications, there could be a service level constraint requiring a high percentage of orders to be delivered on time. This scenario can be modeled using an indicator variable for each product as $I_{p} \in \{0,1\}, \forall p \in P$. Constraints \eqref{ind2} are used to determine if a deadline is violated. Although the GE Gas Power business does not employ a penalty cost or service level approach in meeting deadlines, we illustrate how our model could be modified for these situations.

\begin{subequations} 
\begin{align}
\label{ind1}
&\xi_{p}  - \epsilon_{p} \leq e_{ip} \leq \Xi_{p} + \ell_{p}   , \quad   i \in D, p \in P \\
\label{ind2}
& \ell_{p} \leq \gamma I_{p}   ,\quad   i \in D, p \in P\\
\label{ind3}
&  \sum_{p \in P}  I_{p} \leq \kappa |P|  
\end{align}
\end{subequations}

The high number of constraints and non-zero coefficients in our model stem from the constraint sets used to define flow balance constraints for the in-transit nodes. Preliminary experiments showed that these constraints cause the solver to spend a significant amount of time solving the LP relaxation during the branch and bound search. To enhance the solution speed of the LP relaxation, we define the variables $r_{ipt} \in \{0,1\}$ to represent whether product $p$ stays at in-transit node $i$ at time $t$. We can then replace Constraints \eqref{in-transit-1}, \eqref{in-transit-2}, and \eqref{in-transit-3} with the following constraints:

\begin{subequations}
\begin{align}
\label{invenotry-in-transit-1}
 &   \sum_{ j : (j,i) \in E} \sum_{m \in M^E }     x_{ijmpt} + r_{ipt} = r_{ip(t-1)}  + \sum_{ k: (i,k) \in E} \sum_{m \in M^E } x_{kimp(t - \tau_{kim})}
     \quad  i \in I,  p \in P,    t \in T  \\
%      \label{invenotry-in-transit-2}
% &      r_{ipt} \leq \max_{j \in S} \{s_{jp}\}  \Psi_{ipt}  \quad i \in I,  p \in P, t \in T \\
 \label{invenotry-in-transit-3}
 & \sum_{t \in T}   r_{ipt}  \leq \rho  \quad i \in I,  p \in P
\end{align}
\end{subequations}

Constraints \eqref{invenotry-in-transit-1} are the traditional inventory balance constraints. The LHS indicates that a product either remains at its location at time $t$ as inventory or departs for its next destination. Similarly, the RHS expresses that a product either remains at the location at $t-1$ or is transported from a different location. %Constraints \eqref{invenotry-in-transit-2} decide whether a product stays at the location as inventory, acting as big-M constraints. We use $\max_{j \in S} {s_{jp}}$ as the big-M value.
 Constraints \eqref{invenotry-in-transit-3} ensure that a product does not stay in inventory for more than $\rho$ time periods. According to our experiments with real-world datasets, these constraints perform better than those presented in the IP model.

Lastly, the B$\&$B process may suffer due to symmetry among certain variables, hence, symmetry-breaking constraints may be useful. Under the assumption that FCL units departing from the same port share the same characteristics, including capacity and fixed costs, we can number each container unit in $M^F$ for each port and introduce the following constraints to cope with the symmetry (see Constraints \eqref{symmetryBreaking}). These constraints ensure that if more than one container is to be booked in the same time period following the same route, the model will prefer booking the FCL unit with the smaller number.

\begin{equation} \label{symmetryBreaking}
    z_{ij\text{"FCL''}_kt} \geq  z_{ij\text{"FCL''}_{k+1}t}, \quad   i \in I,  j \in I,   t \in T, k \in \{ 1, \cdots,|M^F|-1\}
\end{equation}

\section{A Knapsack-based Heuristic Approach with a Rolling Horizon} \label{A Knapsack-based Heuristic Approach}

The integer program formulation discussed in the previous section takes a long time to solve real-world problems. Therefore, we propose a novel heuristic algorithm that utilizes the knapsack algorithm to quickly generate a solution. This algorithm operates on a forward-and-back-planning approach, commencing from the current week and determining the actions to take each week for the entire planning horizon. Each week, a decision is made on the number of FCL containers to use and the products that need to be placed in them. Upon reaching the end of the planning horizon, any products that could not be placed in FCL containers are shipped in LCL containers or by air, depending on transportation costs and feasibility.

The algorithm begins by computing the feasible start and end dates for using an FCL container, an LCL container, and air transportation for each product. These dates hinge on the date the product is available to be shipped, and its earliest and latest required delivery dates. For simplicity, we assume that $t_{air}, t_{ocean}$ and $t_{ground}$ represent air, ocean, and ground transportation times in weeks and are the same for each transportation mode across the network. This assumption holds for real-world data when utilizing the optimization model. However, should the transportation times vary, our heuristic would still function with a minor modification. One would need to create individual time ranges for each FCL and LCL unit available in the network. For product $p$, the time range for air transportation can be set as $t_A[p] = ( \max(a_{p}, \xi_{p}- t_{air}),\Xi_{p} - t_{air} )$, taking into account deadlines and availability times. Conversely, the LCL time range is calculated as $t_L[p] = ( \max(a_{p} + t_{ground}, \xi_{pi}- t_{ocean}- t_{ground} - \rho), \Xi_{p} - t_{ocean}- t_{ground})$. Unlike air, we consider here the time taken to reach and leave the port via ground transportation. In addition, the lower bound of the range accounts for the idle wait time at intermediate nodes. Lastly, the FCL time range $t_F[p]$ includes an additional term of $\upsilon$ within the max function to consider the FCL booking time constraint as presented in Constraints \eqref{restrictFCL}.

The algorithm then calculates the lowest alternative cost of shipping each product without placing it in an FCL container. This cost is the minimum of the LCL cost and the air shipment cost. The transportation cost of product $p$ via FCL, LCL, and air are denoted by $F_{m}[p]$, $L[p]$, and $A[p]$, respectively. The sub-index $m$ represents each individual FCL unit connection in $M^F$ (see Table \ref{Table:1}), defined as the edges $(i,j)$ connected via FCL connections (e.g., FCL unit between Shanghai and Greenville Ports). Note that if multiple FCLs are available between nodes $i$ and $j$, we duplicate the existing edge in the network accordingly. To calculate these values, we generate induced networks as $G_{L}^p$, $G_{A}^p$, and $G_{F_{m}}^p$, each associated with the logistics costs of product $p$. The definitions for these networks are provided below:
\begin{enumerate}
\item $G_{L}^p$: An induced network consisting of ground and LCL connections.
\item $G_{A}^p$: An induced network solely consisting of air connections.
\item $G_{F_m}^p$: An induced network consisting of ground and FCL connections, where the associated costs are computed as $c_{ijmp}$ and 0, respectively. Note that fixed costs for FCL units (i.e., $f^F_{m}$) are not included in the cost calculation.
\end{enumerate}

In these definitions, the costs are calculated as the product's cost across each respective network, providing a realistic representation of the logistics costs associated with each product. We solve the shortest path problem using Dijkstra's algorithm, which returns the route yielding the cheapest cost for each induced network. Variables $F_{m}[p]$, $L[p]$, and $A[p]$ are set as the cost of those routes. In total, we solve $(|M^F| + 2)|P|$ shortest path problems as preparation for the Knapsack problems. The value of each $p$ for each FCL unit is set as $v_{m}^p = \lfloor\min(L[p], A[p]) - F_{m}[p]\rfloor$, where the floor function is applied to obtain an integer value. The goal is to examine how much savings could occur when $p$ is placed in an FCL unit rather than being sent via air or LCL. The value $v_{m}^p$ cannot be a negative value unless such a product is pre-processed before starting the solution process.

Once all these calculations are completed, we iterate over the time period from $\upsilon$ to $|T|-t_{ocean}- t_{ground}$ to decide when to book an FCL unit and which products should be placed in the container. At each iteration, we first identify the products that can go with the FCL based on $t_F[p]$ and solve $|M^F|$ knapsack problems. The knapsack problem involves a weight limit/capacity and a set of items, where each item has a weight and a value. The goal is to determine which items should be selected without exceeding the total weight capacity of the knapsack, while maximizing the sum of the values of the selected items. Several solutions to the knapsack problem are available in the literature, including the dynamic programming approach, the branch and bound approach, and hybrid approaches. In our case, we utilize the dynamic programming algorithm to solve the knapsack problem.

Let the objective value of each individual knapsack problem be represented by $o_{m}$, for all $m \in M^F$. We sort each container by $o_{m} - f^F_{m}$ in ascending order and choose the container with the highest non-negative $o_{m} - f^F_{m}$. If the objective value is less than the fixed cost (i.e., $o_{m} - f^F_{m} \leq 0$), it indicates that opening such an FCL unit is not beneficial.

It is important to note that the value provided by the product represents the least cost of the alternative shipment. If the total value obtained from all the products placed in an FCL container using the knapsack algorithm is greater than the fixed cost of shipping the FCL container, then it is worthwhile to book an FCL container for shipping those products. However, there is a possibility of obtaining a higher value by waiting for one more period and shipping the FCL container in the next period. In the next period, there may be more products available for shipping via FCL, and some products may no longer be feasible to ship by FCL due to deadline or capacity constraints.

To evaluate the products available for shipping via FCL in the next period, we use the knapsack algorithm again to compute the value obtained from all the products assigned to the FCL container. We subtract the value of the products that were assigned to the FCL container in the previous period but are no longer eligible for FCL shipment from the total value of products assigned to the FCL container in the current period. If this net value is greater than the value of the FCL products in the previous period, an FCL container will not be shipped in the previous period.

We then determine the value of postponing the shipment of the FCL container to the following period. This process is repeated until the value of shipping the FCL container in a period becomes less than that of the previous period. The algorithm ultimately decides to ship the FCL container in the period with the highest value for the products shipped, minus the value of the products that were placed in an FCL container in earlier periods but are no longer eligible for FCL shipment.

From a mathematical perspective, let us consider the scenario where we are at time $t^{\prime}$ and select the FCL unit $m^{\prime}$ with the highest non-negative value. Before making a final decision, we move to the next time period ($t^{\prime}+1$) and repeat the same process with a slight modification. When calculating each $o_{m}- f^F_{m}$ value at $t^{\prime}+1$, we subtract not only $f^F_{m}$ but also the value of each product that is included in $m^{\prime}$ at $t^{\prime}$ but is excluded at $t^{\prime}+1$ if the corresponding FCL unit is booked. This modification helps us capture two factors: i) how much we lose if certain products are not consolidated, and ii) how much we save if other products are consolidated and become available at a later time period.

If the best value $o_{m}- f^F_{m} - \sum_{p \in m^{\prime}} v_{m^{\prime}}^p$ obtained at $t^{\prime}+1$ is worse than $o_{m^{\prime}}- f^F_{m^{\prime}}$, then we go back to time $t^{\prime}$ and open the FCL unit $m^{\prime}$. If not, we continue marching forward in time and repeat the same process until the total value stops increasing. It's important to note that a product may be in different container units at different time periods as we move forward in time without being excluded from consolidation. Whenever the corresponding product is excluded in a future time period, we use $v_{m^{\prime}}^p$ associated with the candidate container unit selected at time $t^{\prime}$. In the end, we select the container unit that brings the highest return within the evaluated time periods.

If an FCL booking decision is made and there are still products left out, along with more available containers, the algorithm goes back to time $t^{\prime}$. We restart the knapsack solution process to decide if it makes sense to book another container. We also ensure that during the same time period: a) no same FCL unit is booked more than once, and b) no more than $\lambda$ FCL units are booked from the same port.

Once the time horizon is completed, we decide whether a product should be shipped by air or LCL in case it is not placed in a container. This decision is made by examining $\min(L[p], A[p])$ along with $t_A[p]$ and $t_L[p]$. As a result, we obtain a feasible logistics planning decision. In the Online Supplement, we provide the pseudo-code of the algorithm (see Section A). This heuristic offers two notable benefits. First, it does not require any commercial solver support, meaning that using the algorithm does not incur any additional costs for business owners. Second, it can be used to warm start the B$\&$B tree if a solver is utilized, which has been shown to be successful in various optimization studies \citep{egger2021warm, camur2022star}.

It is worth mentioning that if the dynamic programming approach becomes a bottleneck during the solution process, one might consider solving the knapsack problem as an IP model using an open-source solver or utilizing a greedy heuristic that provides a near-optimal solution. In future studies, this heuristic can be further improved by examining multiple Asian lanes (e.g., China $\&$ Singapore to the U.S.) and/or designing approximation guarantees \citep{vogiatzis2019identification, hermans2022exact}.

\section{Computational Experiments} 
\label{Computational Experiments}

We perform  the computational experiments using the Python API and Gurobi solver 9.5.1 on a computer having an $11^{th}$ Gen Intel Core i7-11850H processor and 32 GB of RAM. For the heuristic approach, we utilize the \textit{joblib} library to parallelize both Dijkstra's algorithm and the dynamic programming methods to enhance the solve time. We assure the maintenance of confidentiality by using synthetic test cases. In the rest of this section, we explain the process of creating these test cases and provide a performance comparison between the mathematical formulation and the heuristic algorithm on these test cases. %Our goal is to offer insights and potential strategies for practitioners by testing the algorithms on various realistic scenarios. 

\subsection{Data Set and Scenario Discussion} \label{Data Set 
Discussion}

% To generate our test cases, we employ a combination of high-level histograms, distributions, and averages derived from historical data, along with certain assumptions. The aim is to replicate the statistical variations observed in real-world scenarios. These settings enable us to 

We create randomized data across dimensions such as product weight, volume, origin-destination pairs, ready-to-ship dates, earliest and latest requirement dates, number of product orders, and the length of the planning horizon. We describe the generation process for weight and volume data below.

Gross Weight: For 63$\%$ of the shipments, the weights are uniformly distributed between 50 and 500 kilograms. The remaining 37$\%$of shipments have weights uniformly distributed between 500 and 5,000 kilograms.

Volume: The volume of a shipment is determined by multiplying its gross weight in metric tons by a coefficient. The coefficient is assigned as follows:

\begin{itemize}
    \item 		For 12.5$\%$ of the products, the coefficient is uniformly distributed between 0.5 and 1.5.
    \item 	For the next 67.5$\%$ of the products,  the coefficient is uniformly distributed between 1.5 and 3.0.
    \item		For the next 10$\%$ of the products, the coefficient is uniformly distributed between 3.0 and 4.5.
    \item 	For the next 2.5$\%$, the coefficient is uniformly distributed between 4.5 and 6.
    \item For the remaining 7.5$\%$ of the products, the coefficient is uniformly distributed between 6 and 7.5.
\end{itemize}

Air-Charge Weight: The air-charge weight is determined based on the volume-to-gross-weight ratio of the shipment. If the ratio is less than 3, the air-charge weight is set equal to the gross weight of the shipment. Otherwise, it is calculated as 1.2121 times the gross weight.

Origin-Destination Pairs: The origin of all shipments is set as an inland location in China. The number of final destinations varies between 1, 2, and 3. These destinations are selected from the cities of Greenville, Bangor, Atlanta, and Schenectady, where GE Gas Power has manufacturing plants and/or warehouses. For exporting parts by ocean, the ports of Shanghai and Qingdao in China are used, while for receiving them in the US, the ports of Baltimore, Charleston, Newark, and Savannah are utilized.

The destination of each shipment is randomly determined from the available options based on the specific case being studied. In cases where multiple alternative destinations are available, each destination has an equal likelihood of being chosen.

Ready-to-ship Dates and Deadlines: The week when a product becomes available to be shipped is uniformly distributed over the first 35$\%$ of the planning horizon. The earliest acceptable need-by-week is determined based on the ready-to-ship week using the following rules:

\begin{itemize}
\item	There is a 7.5$\%$ chance that the shipment will be needed 11 weeks after its ready-to-ship date.
\item	There is a 40$\%$ chance that it will be needed 12 weeks after its ready-to-ship date.
\item	Otherwise, the time difference between the ready-to-ship date and the earliest need-by date is randomly distributed between 13 weeks and the end of the optimization planning horizon minus 2 weeks.
\end{itemize}

The latest acceptable need-by-week is set as the earliest need-by-week plus 2 weeks. 

Costs: Booking an FCL unit incurs a cost of $\$12,538$, while an LCL unit has a bunker cost of $\$1,160$. Additionally, there is a linear cost per cubic meter (CBM) for LCL units, which is uniformly distributed between $\$700$ and $\$900$. For air and ground transportation, the costs are $\$13.23$ and a random cost between $\$0.1$ and $\$0.5$ per kilogram (kg), respectively. For example, if a product weighs 10 kg and is transported by air, the total cost would be $\$132.3$.

Transportation Times: Both ground and air transportation are assumed to take two weeks, while ocean transportation takes seven weeks.

Capacity: The capacity of an FCL container is set to 20 tons. LCL and air shipments are not subject to capacity constraints.

Number of Products: The number of purchase order (PO) line items per test case varies among the following options: $\{$50, 100, 250, 500, 1000$\}$.

Planning Horizon: The planning horizon is chosen to be either 6 or 12 months long.

The test cases generated are numbered as presented in Section B of the Online Supplement. We solve each test case under three different scenarios:

\begin{itemize}
    \item[a)]	Baseline Scenario: In this scenario, the model is allowed to utilize both ports in China (Shanghai and Qingdao) as well as the four ports in the U.S. (Baltimore, Charleston, Newark, and Savannah).
    \item[b)] 	Port Closure Scenario: This scenario examines the impact of the closure of the Shanghai Port, which occurred in early 2022. 
    %We analyze how the logistics planning decisions are affected by the unavailability of this port.
    \item[c)] 	NO FCL utilization Scenario:  In this scenario,  the use of FCL containers is prohibited.
\end{itemize}

\subsection{Warm Start Analysis}
In order to evaluate the impact of supplying the solution returned by the heuristic algorithm as an initial solution to the Gurobi solver, we conduct experiments using a time limit of one hour and an optimality gap limit of 0.5$\%$. We solve each of the 30 test cases both with and without warm start (WS) using Gurobi. The optimality gaps and solution times were compared to assess the effectiveness of the WS strategy (see Section C in the Online Supplement for details).

Based on our analysis, we conclude that warm starting the Gurobi solver with the solution obtained from the heuristic algorithm did not significantly improve the solution quality or solution time in most instances. The performance was generally similar between the two approaches, with a few exceptions. For instance, in test case nineteen, the optimality gap increased from 44.6$\%$ to 76.1$\%$ when using the WS strategy. This behavior is somewhat expected, as WS is typically beneficial when the solver struggles to find initial feasible solutions. Therefore, for the remainder of the computational experiments, we conducted them without utilizing the WS strategy.

\subsection{Computational Results}
\label{Computational Results}
In this section, we compare the results of the IP solution with the Knapsack-based heuristic using the 30 test cases under  the three scenarios  described earlier. All the computational results are presented in Table~\ref{table1}.  For each test case and  scenario, we  report the following information: a) Minimum total cost found by the IP solution, b) optimality gap of the IP solution, c) heuristic error, which represents the percentage increase in the objective function value compared to the IP solution, and d) duration of the run times for the IP and heuristic solutions.

For the "Port Closure" and "No FCL" scenarios, we also report the percentage cost increase compared to the baseline scenario. The heuristic error is computed as $\frac{(obj_{\text{Heur.}}-obj_{\text{IP}})}{obj_{\text{IP}}}$, where $obj_{\text{IP}}$ and $obj_{\text{Heur.}}$ denote the objective function values returned by the optimization model and the heuristic algorithm, respectively.

Additionally, we set a time limit of one hour for the Gurobi Solver when solving the IP. If the time limit is reached before reaching the optimal solution, we report the best feasible objective function value found along with its optimality gap. We indicate in the table with "TL" to signify that the time limit was reached.

\begin{sidewaystable}
%\begin{table}[!hp]
  \begin{center}
\setlength\tabcolsep{0pt}
\renewcommand{\cellset}{\bfseries\linespread{1}\selectfont}
\renewcommand\cellalign{bc}
  \centering
  \caption{The computational results}{
 \small
 %\begin{adjustbox}{angle=90}
    \begin{tabular*}{\linewidth}{@{\extracolsep{\fill}} 
                        *{1}{c}
                        *{1}{r r r r r}
                        *{2}{r r r r r r}}
     \toprule
     
    & \multicolumn{5}{c}{Baseline}   
    & \multicolumn{6}{c}{Port Closure (Covid Impact)}           
    & \multicolumn{6}{c}{No FCL utilization }     \\
    
    \cmidrule(l){2-6}
    \cmidrule(l){7-12}
    \cmidrule(l){13-18}
    
    {\makecell{Test \\ Case}}
    & {\makecell{IP \\ Obj. \\ $(\$)$ }} & {\makecell{Opt.\\ Gap \\ $(\%)$}} & {\makecell{Heur.  \\Error \\ $(\%)$}}   
    & {\makecell{IP\\ Time \\ (sec)}} & {\makecell{Heur.\\ Time \\ (sec)}}  & {\makecell{IP \\ Obj.  \\ $(\$)$}} & {\makecell{Cost \\ Inc. \\ $(\%)$ }} & {\makecell{Opt.\\ Gap \\ $(\%)$}} & {\makecell{Heur.  \\Error \\ $(\%)$}} 
    & {\makecell{IP\\ Time \\ (sec)}} & {\makecell{Heur.\\ Time \\ (sec)}} &
     {\makecell{IP \\ Obj. \\ $(\%)$}} & {\makecell{Cost\\ Incr. \\ $(\%)$}} & {\makecell{Opt.\\ Gap \\ $(\%)$}} & {\makecell{Heur.  \\Error \\ $(\%)$}}   
    & {\makecell{IP\\ Time \\ (sec)}} & {\makecell{Heur.\\ Time \\ (sec)}}\\
    \midrule
1&	54328.34&	0.0&	11.3&	47.2&	8.5	&81540.76&	50.1&	0.0&	6.9&	55.7&	6.6&	192792.55&	254.9&	0&	0&	1.0	& 2.5\\
2&56,879.57&0.0&19.8&9.05&5.08&85,251.32&49.9&0.0&12.4&5.33&3.44&198,860.85&249.6&0&0&0.8&0.52\\
3&42,238.66&0.0&18.4&6.58&5.52&60,150.91&42.4&0.0&13.0&1.46&4.43&166,685.02&294.6&0&0&0.79&0.49\\
4&82,383.24&0.0&4.2&5.69&6.02&102,191.90&24.0&0.0&3.5&1.68&4.1&167,233.24&103.0&0&0&1.39&0.5\\
5&84,828.06&0.0&14.1&3.87&5.68&105,676.47&24.6&0.0&11.3&1.55&3.61&196,405.07&131.5&0&0&1.38&0.5\\
6&93,086.68&0.0&1.5&4.1&6.11&114,180.22&22.7&0.0&1.2&1.38&3.81&175,369.93&88.4&0&0&1.27&0.55\\
7&77,421.90&3.2&7.9&TL&10.7&118,584.99&53.2&0.6&5.3&TL&7.12&351,666.21&354.2&0&0&2.11&0.84\\
8&86,369.39&7.2&20.9&TL&9.64&128,195.82&48.4&3.5&2.0&TL&6.83&340,399.16&294.1&0&0&2.18&0.77\\
9&79,488.39&2.2&4.6&TL&10.08&119,197.48&50.0&0.5&4.3&1686.4&6.47&321,042.91&303.9&0&0&1.82&0.71\\
10&106,144.57&0.0&4.8&6.91&10.72&144,344.69&36.0&0.0&3.6&4.86&6.41&327,754.27&208.8&0&0&3.55&0.82\\
11&129,833.22&0.5&12.2&61.62&10.96&184,255.65&41.9&0.5&8.5&26.1&6.93&387,683.99&198.6&0&0&3.52&0.99\\
12&120,856.16&0.4&13.5&18.42&10.57&167,224.38&38.4&0.5&9.7&8.21&6.35&366,863.21&203.6&0&0&3.12&0.98\\
13&218,513.62&1.9&5.7&TL&29.29&339,695.00&55.5&1.1&1.2&TL&18.43&907,054.77&315.1&0&0&6.03&2.46\\
14&187,894.68&2.3&2.5&TL&28.1&286,579.95&52.5&1.5&4.5&TL&16.58&753,233.42&300.9&0&0&5.32&2.11\\
15&220,506.30&1.0&8.7&TL&28.69&337,706.01&53.2&0.7&5.8&TL&18.68&881,967.26&300.0&0&0&5.14&1.94\\
16&205,525.33&2.0&15.8&TL&24.7&319,261.86&55.3&1.3&8.8&TL&16.53&877,975.83&327.2&0&0&9.21&2.36\\
17&215,729.25&2.8&12.5&TL&25.18&326,375.35&51.3&1.6&8.2&TL&20.5&871,101.59&303.8&0&0&8.34&2.27\\
18&216,213.11&3.2&15.8&TL&25.51&328,760.93&52.1&2.1&8.5&TL&18.05&915,945.01&323.6&0&0&6.79&1.86\\
19&451,192.02&44.6&11.7&TL&76.09&712,970.74&58.0&0.4&9.4&1905.08&52.51&1,811,405.86&301.5&0&0&12.17&3.54\\
20&498,615.52&74.3&9.9&TL&78.73&787,386.86&57.9&52.0&9.7&TL&52.5&1,865,309.00&274.1&0&0&11.98&3.73\\
21&521,944.95&76.1&2.1&TL&73.35&752,886.74&44.2&49.2&11.5&TL&40.49&1,823,922.08&249.4&0&0&10.42&3.05\\
22&397,965.61&48.5&5.1&TL&54.31&621,873.89&56.3&0.5&2.4&2075.74&25.9&1,743,223.80&338.0&0&0&21.78&3.72\\
23&412,912.64&47.0&5.5&TL&57.11&630,800.07&52.8&0.5&5.3&3420.82&29.48&1,710,630.75&314.3&0&0&20.16&3.94\\
24&423,364.25&48.9&9.6&TL&58.1&650,735.61&53.7&0.5&5.1&2300.2&26.79&1,787,195.80&322.1&0&0&16.92&2.89\\
25&1,235,003.22&85.8&3.3&TL&237.69&2,145,037.72&73.7&67.3&1.5&TL&87.9&3,641,673.11&194.9&0&0&29.51&6.17\\
26&1,157,173.31&81.0&7.6&TL&223.62&2,069,246.17&78.8&67.0&2.9&TL&80.04&3,452,146.96&198.3&0&0&26.11&5.52\\
27&1,373,587.07&83.0&1.1&TL&192.97&2,312,049.14&68.3&69.4&2.3&TL&85.96&3,653,723.88&166.0&0&0&28.87&5.42\\
28&839,723.93&72.1&9.0&TL&111.05&1,314,566.58&56.5&7.0&7.9&TL&59.98&3,684,769.25&338.8&0&0&58.85&5.7\\
29&861,949.79&67.2&6.2&TL&106.42&1,298,403.38&50.6&39.8&4.9&TL&55.93&3,453,393.49&300.6&0&0&49.04&4.71\\
30&868,716.92&67.2&7.6&TL&103&1,307,662.40&50.5&41.2&6.2&TL&55.7&3,609,491.22&315.5&0&0&44.3&4.09\\

\textbf{Average} & & \textbf{27,4} & \textbf{9.1}&\textbf{ 2524.4} &\textbf{54.4} & &\textbf{50.1} & \textbf{13.6}& \textbf{6.3} & \textbf{2303.2} & \textbf{27.6} & &\textbf{262.3} & \textbf{0} & \textbf{0} &\textbf{13.1}& \textbf{2.5} \\

    \bottomrule
    
    \end{tabular*}
} \label{table1}
   \end{center}
%\end{table}

\end{sidewaystable}

In the baseline scenario, the solver performs satisfactorily within the time limit for all test cases with 50, 100, and 250 products (test cases 1-18). It either reaches the optimal solution or achieves an optimality gap below the limit (0.5$\%$) in almost 40$\%$ of these instances. %For the remaining cases in this group, the average optimality gap is 2.4$\%$.
However, for test cases with 500 or more products, the average optimality gap increases to 66.3$\%$. As expected, the solver's performance decreases as the problem size increases.

When comparing the two approaches, the IP algorithm with a one-hour time limit outperforms the heuristic in terms of solution quality. On average, the heuristic returns recommendations with an additional total cost of 9.1$\%$. Although the IP algorithm outperforms the heuristic in every test case in terms of the solution objective, the heuristic error decreases as the problem scale increases. For test cases where the solver consumes the entire time allowed and the optimality gap exceeds 40$\%$, the average heuristic error drops to approximately 6$\%$. 

%This suggests that as the problem becomes larger, the heuristic may start to approach the performance of the commercial solver.

In terms of solution times, the heuristic algorithm completes in less than four minutes, with an average runtime of 54.4 seconds across all test cases. In comparison, the solver fails to complete within the one-hour time limit for 70$\%$ of the test cases. It is important to note that even though there may be a large gap between the best feasible solution found and the theoretical lower bound, it does not necessarily imply poor solution quality. It is possible that an optimal or near-optimal solution is obtained, but there was insufficient time for the solver to prove optimality and improve the primal bounds.

\begin{figure}[htbp]
    \begin{minipage}{0.48\textwidth}
   \includegraphics[page=1,width=1\linewidth]{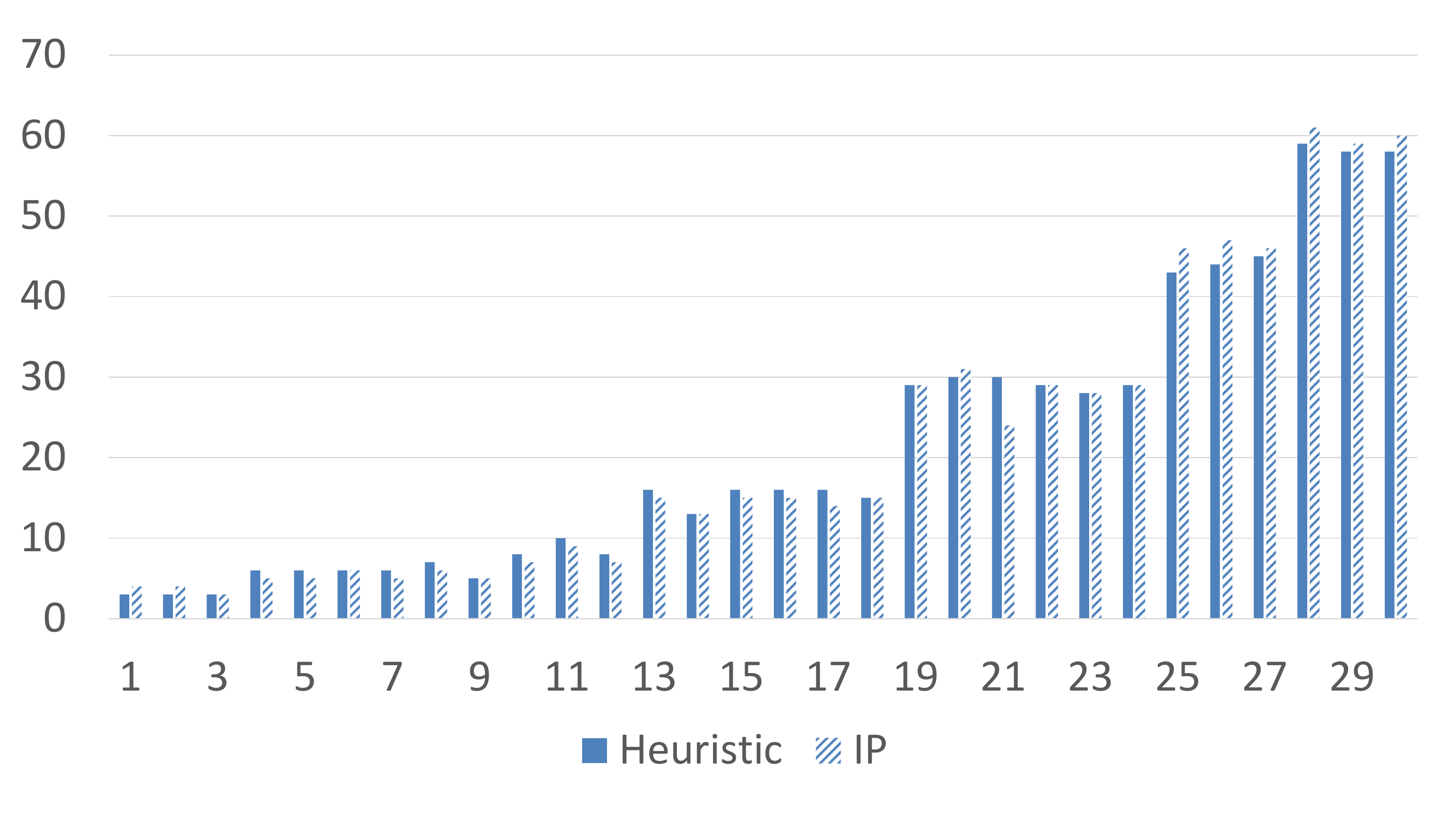}  
\caption{Comparison of number of FCL units booked  by the heuristic and IP algorithms in the baseline scenario}
  \label{Fig3}
    \end{minipage}%
      \hfill
    \begin{minipage}{0.48\textwidth}
   \includegraphics[page=2, width=1\linewidth]{Fig.pdf}  
\caption{Comparison of the ground transportation costs returned by the heuristic and IP algorithms in the baseline scenario}
  \label{Fig4}
    \end{minipage}
\end{figure}

As we discusses in Section \ref{Introduction}, our primary goal is to identify the right FCL booking decisions within a short amount of time rather than competing with commercial solvers with respect to the solution quality. Thus, we compare the number of FCL units booked in each solution within the baseline scenario. The heuristic and IP model book the same number of containers in nine out of thirty instances (see Fig.~\ref{Fig3}). The largest difference is observed in test case 21, where the heuristic recommends six additional containers. However, the error in the objective is reported as only 2.1$\%$. This discrepancy arises because the IP model sends 35 products using LCL units, resulting in a total cost of $\$$156,660.37, while the heuristic only sends two products via LCL. Due to the different ports used, the costs of ground transportation vary, leading to similar objective values.

Overall, our heuristic algorithm makes the right booking decisions in terms of timing in most cases. The cost difference primarily arises from the products placed in each container, which in turn impacts the ground transportation costs based on the final destinations. In fact, in the solutions obtained by the heuristic, the ground transportation costs are higher in 80$\%$ of the test instances (see Fig.\ref{Fig4}). For the test cases where the IP model has a higher ground transportation cost, the maximum cost difference is approximately $\$$23,000.

\begin{figure}[htbp]
    \begin{minipage}{0.48\textwidth}
   \includegraphics[page=3,width=1\linewidth]{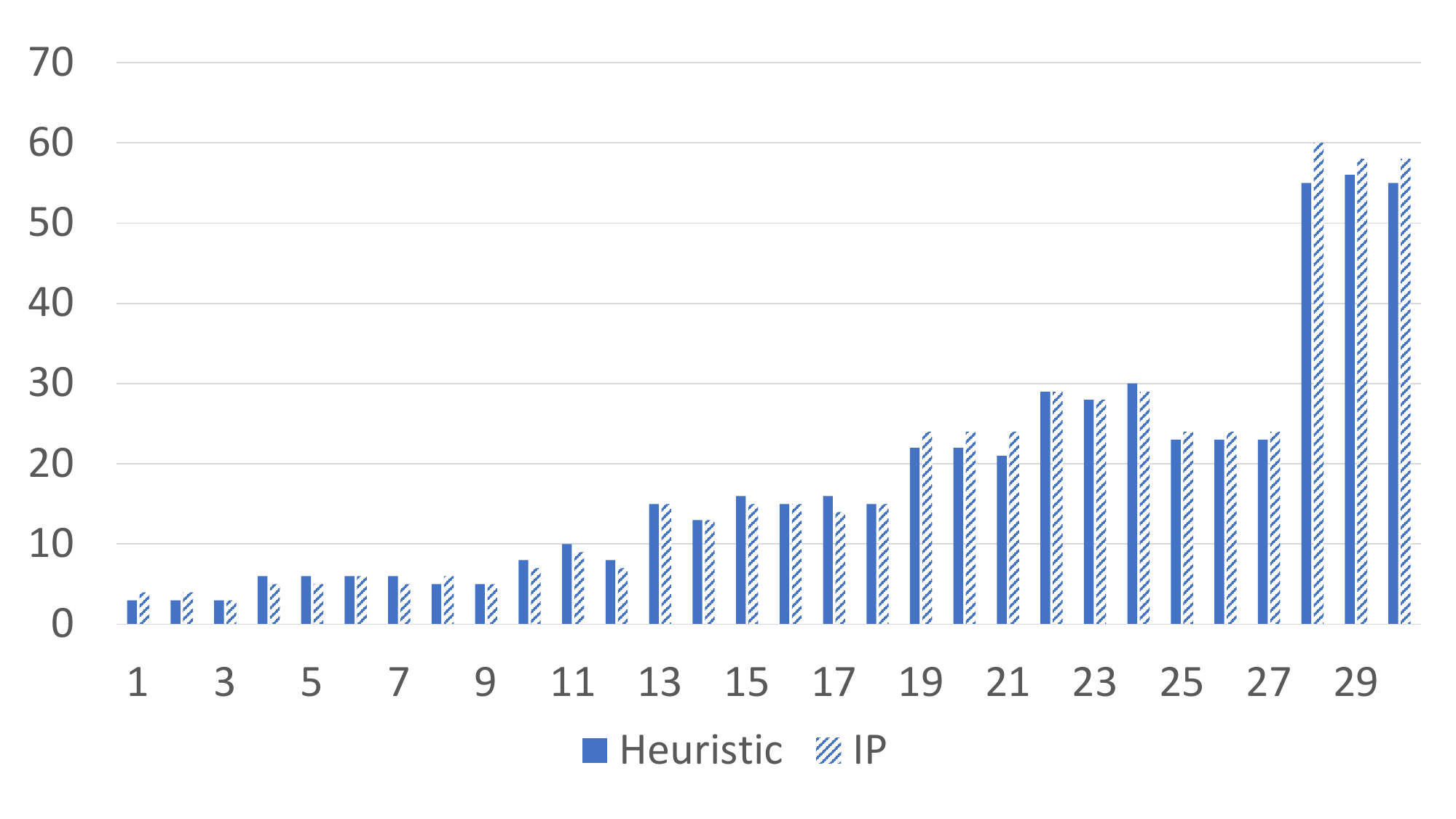}  
\caption{Comparison of number of FCL units booked by the heuristic and IP algorithms in the ``Port Closure" scenario}
  \label{Fig5}
    \end{minipage}%
      \hfill
    \begin{minipage}{0.48\textwidth}
   \includegraphics[page=4, width=1\linewidth]{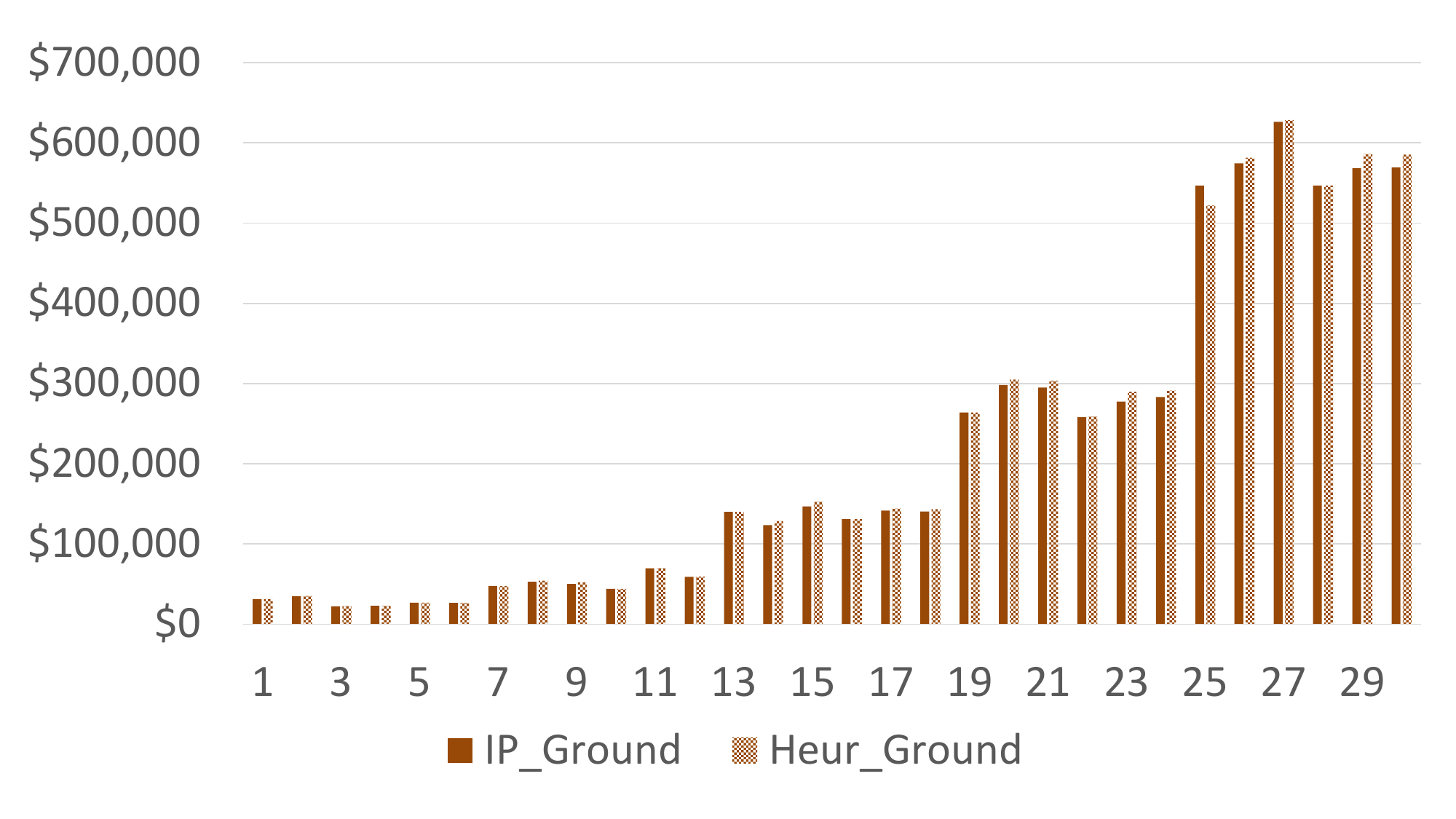}  
\caption{Comparison of the ground transportation costs returned by the heuristic and IP algorithms in the ``Port Closure" scenario}
  \label{Fig6}
    \end{minipage}
\end{figure}

 In the "Port Closure" scenario, similar to the baseline scenario, the IP approach outperforms the heuristic in terms of the objective value, with an average heuristic error of 6.3$\%$. Since there are fewer opportunities for FCL bookings due to the port closure, the heuristic demonstrates better performance compared to the baseline scenario.  As shown in Fig.~\ref{Fig5}, a smaller number of FCL units are utilized in this scenario. Additionally, the standard deviation of the difference in the number of FCL units booked is 1.49, which is better than the baseline scenario with a standard deviation of 1.62. Furthermore, since similar ports are used, the ground transportation costs are not significantly different, with an average cost difference of approximately $\$2,524$. This represents a 77$\%$ reduction compared to the baseline scenario (see Fig.~\ref{Fig6}). 

\begin{figure}[htbp]
    \begin{minipage}{0.48\textwidth}
   \includegraphics[page=5,width=1\linewidth]{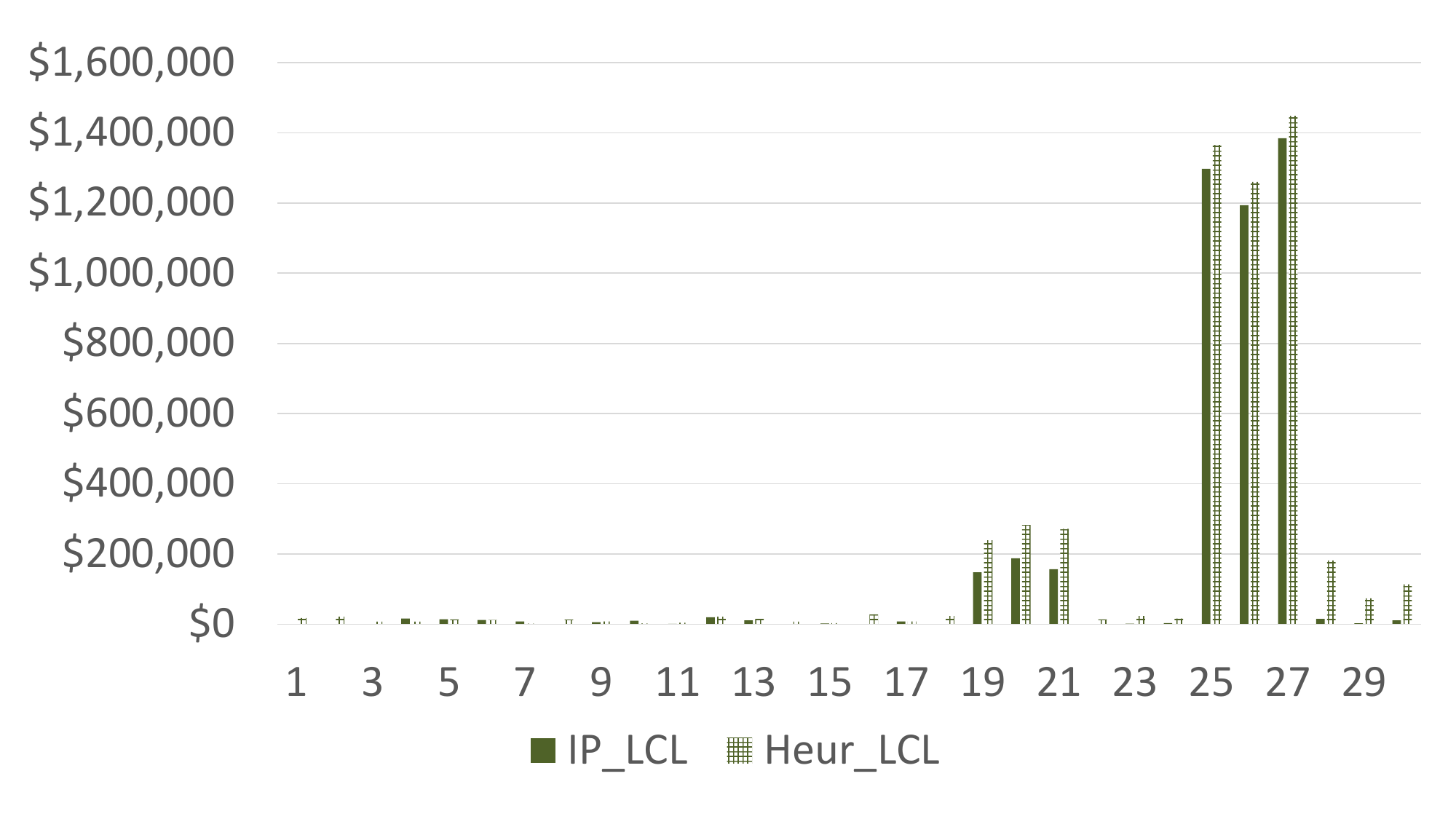}  
\caption{Comparison of the LCL costs  returned by the heuristic and IP algorithms in the ``Port Closure" scenario}
  \label{Fig7}
    \end{minipage}%
      \hfill
    \begin{minipage}{0.48\textwidth}
   \includegraphics[page=6, width=1\linewidth]{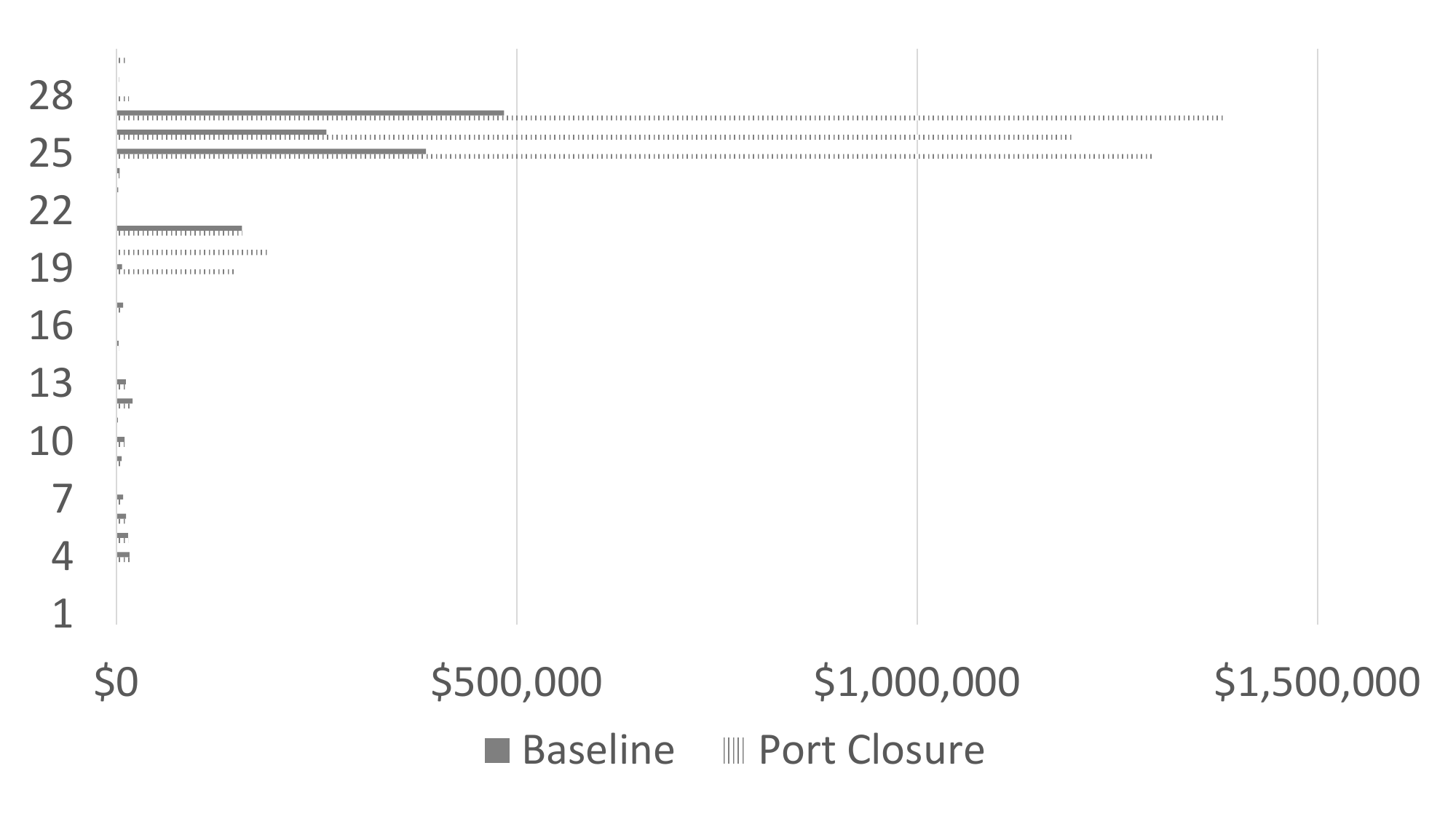}  
\caption{Comparison of the air transportation costs obtained by the IP algorithm in the baseline and ``Port Closure" scenarios}
  \label{Fig8}
    \end{minipage}
\end{figure}

Another factor contributing to the difference in objective values is the handling of products that are not placed in an FCL unit. The heuristic algorithm progresses in time in a rolling manner, and there is a possibility that certain products might not be consolidated into an FCL unit based on the solution returned by the knapsack problem. As a result, consolidation opportunities for those products may be missed. These products are then sent via LCL units or air transport based on the deadline constraints. Consequently, the heuristic incurs higher LCL costs compared to the IP in twenty-six test cases. This leads to an average of $\$$33,283 extra LCL logistics cost per test case, as shown in Fig.~\ref{Fig7}. 

It is important to note that the average total logistics costs increase by $50\%$ compared to the baseline scenario. This is primarily due to the decrease in the number of booked FCL units, resulting in lower FCL costs, which average around $\$35,000$. However, as shown in Fig.~\ref{Fig8}, the total LCL costs exponentially increase in a few test cases. On average, the total LCL costs are approximately $\$103,420$ higher in the "Port Closure" scenario. It is worth mentioning that there is no significant cost difference between the two scenarios in terms of air transportation, which is expected as air transportation incurs high costs compared to ocean transportation options.

As anticipated, due to the smaller number of arcs in the network, the heuristic algorithm runs to completion faster, taking less than one and a half minutes for all cases. On the other hand, the solver uses the entire hour for sixteen out of the thirty test cases. This further highlights the computational efficiency of the heuristic implementation based on our test cases. Particularly, the heuristic can be useful for obtaining initial, high-quality feasible solutions more quickly in cases where the solver struggles to identify such solutions. In fact, there are real-world, large-scale cases where the solver faces difficulties in solving the LP relaxation of the model to initiate the B$\&$B process.

Lastly, in the "No FCL" scenario, the problem becomes much simpler as there is no option for shipment consolidation and prebooking of containers. The decision-maker only needs to choose between air transportation and LCL units for product transportation. In this scenario, both the IP and heuristic approaches provide optimal solutions in every test case without exception. The solver has an average runtime of nineteen seconds and never exceeds one minute. As for the heuristic, it completes on average in 3.3 seconds. This indicates that there is no need for an optimization tool when shipment consolidation does not play a role in the planning process.

\begin{figure}[htbp]
    \begin{minipage}{0.48\textwidth}
   \includegraphics[page=7,width=1\linewidth]{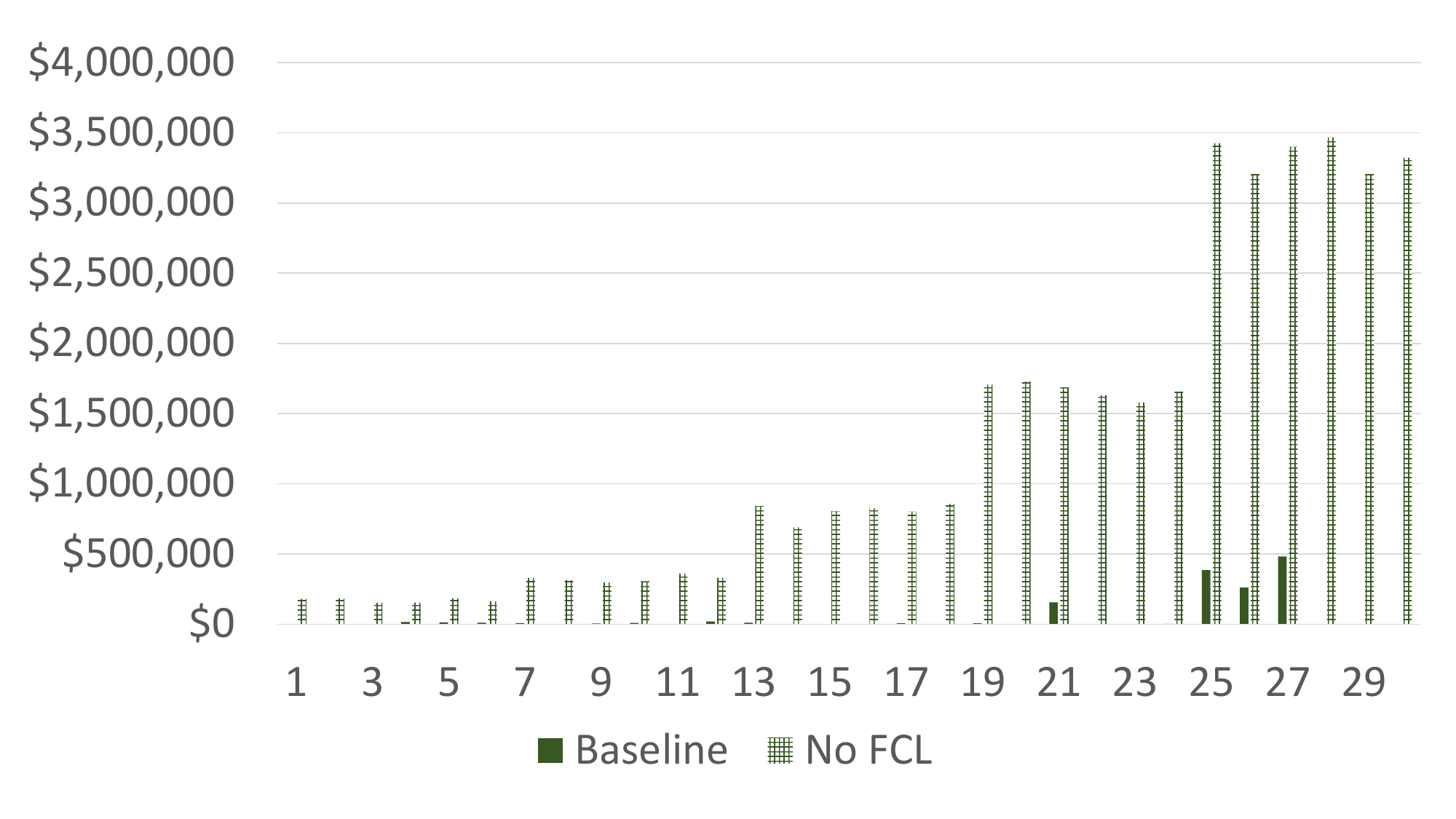}  
\caption{Comparison of the LCL costs  obtained by the IP algorithm in the baseline and ``No FCL" scenarios}
  \label{Fig9}
    \end{minipage}%
      \hfill
    \begin{minipage}{0.48\textwidth}
   \includegraphics[page=8, width=1\linewidth]{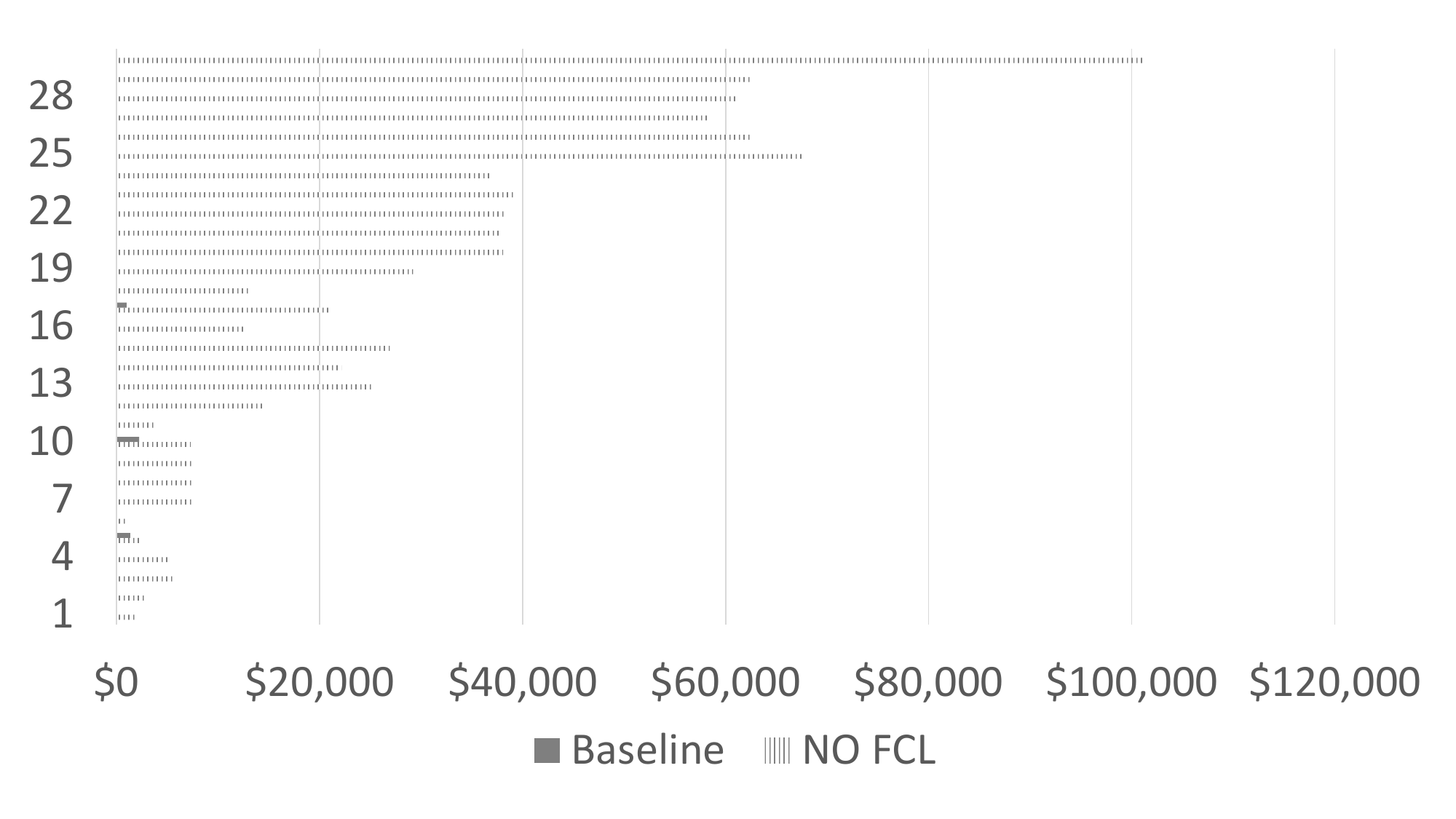}  
\caption{Comparison of the air transportation costs obtained by the IP algorithm in the baseline and ``No FCL" scenarios}
  \label{Fig10}
    \end{minipage}
\end{figure}

The average cost increase between the baseline and "No FCL" scenario is 262$\%$, meaning that the total costs are nearly tripled. This significant cost increase is expected since all the products are shipped via air and LCL, which are more expensive means of transportation compared to a relatively full FCL container. This highlights the importance of FCL consolidations. Shipment consolidation plays a crucial role in reducing total costs in most logistics operations.

For instance, in the baseline scenario, there are nine test cases where no LCL unit is utilized, and air transportation is only used in three test cases. As a result, the average LCL costs increase from $\$$47,000 to $\$$1,260,000, as shown Fig.~\ref{Fig9}. Similarly, when no FCL unit is available, we observe air transportation being utilized in every test case (see Fig.~\ref{Fig10} ) resulting in a 17,876.3$\%$ increase in the average total air transportation costs. These findings demonstrate the cost-saving benefits and efficiency gained through FCL consolidations.

\subsection{Technical and Managerial Insights}
Based on the results of the computational experiments, we can provide some recommendations for different problem sizes and planning scenarios. When the problem scale is small, especially in terms of the number of products, it is recommended to use the IP formulation and solve it with a commercial solver. The IP approach is effective in obtaining optimal solutions for smaller  instances.

As the problem size increases, especially with around a thousand products and larger time horizons, the heuristic algorithm can provide a good solution quickly and is a viable alternative to the IP approach. It can be used as a standalone solution or as a warm start for the commercial solver. The heuristic algorithm demonstrates computational efficiency and can offer high-quality feasible solutions in terms of the FCL booking decisions.

For very large problem sizes with several thousands of products and a complex supply chain network, the proposed heuristic algorithm may outperform the solver and provide a solution in a relatively short amount of time. In these cases, the heuristic can be a valuable tool for obtaining initial feasible solutions efficiently.

Based on the cost analysis conducted in the previous section, it is evident that shipment consolidation offers significant cost-saving opportunities for supply chain and logistics planners. By consolidating multiple products into a single container, both total logistics costs and carbon emissions can be reduced. Therefore, companies should consider incorporating shipment consolidation strategies into their logistics planning processes. Furthermore, considering air consolidation opportunities, similar to ocean container consolidations, may lead to additional cost savings and operational efficiencies.

It is important to acknowledge that careful logistics planning is crucial for ocean transportation, given the considerable time required for shipments. FCL booking decisions are influenced by two critical parameters: deadlines and product availability dates. While deadlines are determined by the problem owner, availability dates come from the suppliers. Considering the uncertainty and disruptions in the global supply chain, future research on predicting product availability dates, as an alternative to relying solely on supplier-provided dates, could greatly enhance the dependability of optimization results.

The benefits and insights obtained from our optimization model can be categorized into strategic, tactical, and simulation aspects. Strategically, the model provides a quarterly overview of the expected assets needed, offering a high-level understanding of resource requirements. Tactically, it can be used on a weekly basis to monitor any changes or deviations from the standard cadence. Finally, the model enables strategic simulations of different scenarios for supply chain planning, allowing decision-makers to evaluate various what-if scenarios and make informed decisions.

\section{Conclusion}
\label{Conclusion}
In this study, we address the inbound logistics and shipment consolidation optimization problem faced by GE Gas Power, focusing on the transportation of gas turbine parts from China to the U.S. We considered various transportation modes, including ground, air, LCL, and FCL, as well as early and late deadlines for the products. Our goal was to develop effective planning strategies for shipment timing and consolidation decisions.

To tackle the problem, we first formulated an integer programming model that incorporated relevant business constraints. We then proposed a rolling horizon heuristic algorithm that leveraged shortest path and knapsack algorithms to provide near-optimal solutions within a reasonable time frame. Computational experiments were conducted to evaluate the performance of the proposed approaches under different realistic scenarios.

The results of our analysis demonstrate the significant cost savings achieved through FCL consolidation, with an average reduction of logistics operations costs by 260$\%$. Having access to multiple ports in both the source and destination locations allows for alternative routes, further reducing costs. The use of alternative container routes also helps in mitigating potential disruptions in the supply chain.

Moving forward, it would be interesting to extend the study to incorporate multiple Asian countries as source locations, considering internal consolidation operations before shipping products to the U.S. However, such practices require careful consideration of border regulations, tariffs, and tax implications. Additionally, optimizing the container load of FCL units could be integrated into the optimization model to make more realistic decisions. Exploring other ocean transportation modes, such as roll-on/roll-off ships, could also be worthwhile to assess their impact on total logistics costs. Researchers may also consider investigating ground and air consolidation opportunities.

%Furthermore, there is room for further improvement in the heuristic approach to enhance the quality of the solutions returned. Future research efforts could focus on refining the heuristic algorithm to provide even higher-quality solutions and explore additional optimization strategies.
% \if0\blind{
% \section*{Acknowledgements}
% The authors acknowledge the generous support from the funding agency of XYZ.	} \fi

\end{document}